\documentclass[]{article}

 \usepackage[]{graphicx}
 \usepackage{epsfig}
 \usepackage{amsmath, amstext, amsfonts}
 \usepackage{amssymb}
 \usepackage{multicol}
 \usepackage[latin1]{inputenc}

\usepackage{anysize} 
\marginsize{2cm}{2.5cm}{1.5cm}{2cm}

\begin{document}

 \title{Scattering Theory for Hamiltonians Periodic in Time}

 \author{J.S. Howland\footnote{University of Virginia, Charlottesville, VA 22903}
  }

 \date{Indiana J. Math. 28, 471-494 (1979)\\
 -- $\cdot$ --\\
{\em \small Received May 1, 1978; Revised September 18, 1978.}
 }



\def\ct{\centerline}
\def\1{\'{\i}}
\def\={\!=\!}
\def\be{\begin{equation}}
\def\ee{\end{equation}}

\maketitle

%

%


{\bf Introduction}. $\quad$ This paper concerns the scattering theory of the Schroedinger evolution equation
\begin{equation}
i\frac{d \psi}{d t}=H(t)\psi~, \qquad -\infty < t < \infty
\end{equation}
where $H$ is a self-adjoint Hamiltonian operator on a Hilbert space ${\cal H}$.

If $U(t,s)$ and $U_0(t,s)$ are the corresponding unitary propagators
that is, the operators that give the solutions of the equations according to the formula
$$
\psi(t)=U(t,s)\psi(s)
$$
then the question reduces to the study of the wave operators
\begin{equation}
W_{\pm}(s)=s-\lim_{t\rightarrow \pm \infty} U_0(t,s)U(t,s)P\, 
\end{equation}
where $P$ is some projection whose introduction may be necessary in order that the limit exist. The first things to be considered
are usually existence of the limits, and the ranges of $W_{\pm}(s)$ and its adjoint. A substantial literature is devoted to the case where $H(t)$ is independent of $t$.

\medskip

In a previous paper, the author considered the case where the difference $H-H_0$ goes to zero suitably as $t\rightarrow \pm \infty$. The method employed was suggested by a procedure of Classical Mechanics for reducing $t$ to a spatial variable. One considers the operator
$$
K=-i\frac{d}{d t}+H(t)
$$
on the space ${\cal H}=L_2(-\infty,\infty; {\cal H})$, and the similar operator $K_0$ corresponding to $H_0(t)$. The study of $W_{\pm}(s)$ was shown to be essentially equivalent to the study of
$$
\Omega_{\pm}(K_0,K)=s-\lim_{\sigma\rightarrow \pm \infty} e^{i\sigma K_0}e^{-i\sigma K} 
$$
In this case, one can take $P=I$. The techniques developed for the time-indep. case, in particular the ``stationary theory'', are then available here.
Actually, we considered in [9], the class of s-adjoint operators $K$ on ${\cal H}$ such that for all smooth scalar functions $\phi(t)$
$$
KM(\phi)-M(\phi)K=-iM(d\phi/dt)
$$
where $M(\phi)f(t)=\phi(t)f(t)$. This symmetry equation asserts {\em formally} that the difference of $K$ and $-id/dt$ commutes with every scalar multiplication, and is therefore an operator-valued multiplication $H(t)$. There exists a propagator $U(t,s)$ for such a $K$, but it may only be a measurable function of $t$ and $s$. Nevertheless, a satisfactory scattering theory (ST), with a unitary $S$-matrix, was obtained. Using the ST, it was easy to show that $U(t,s)$ is strongly continuous for the operators considered in [9]. However, without extra conditions to permit the use of some standard theorem on evolution eqs, it is not at all clear that $U(t,0)\psi_0$ is a strong solution of (1).

\medskip

In any case, questions about strong solutions of evolution equations and the construction of propagators as multiplicative integrals, however interesting in themselves, are largely irrelevant to scattering theory.

\bigskip

In the present paper, we consider the case in which $H(t)$ and $H_0(t)$ are periodic with the same period $a$. This is equivalent to considering operators $K$ and $K_0$ which commute with the unitary time-translation
$$
T_af(t)=f(t-a)~.
$$
Regarding $t$ as a spatial variable suggests decomposing $K$ in the spectral representation of $T_a$, as is done in spatially periodic pbs.
This leads to an abstract correspondence of the scattering problem for (1) to the ST of a certain operator $K$.

\begin{equation}
\tilde{K}=-i\frac{d}{dt}+H(t)
\end{equation}
with periodic BCs on $[0,a]$, or by a Fourier series expansion, as
\begin{equation}
\tilde{K}x_n=2\pi a^{-1}nx_n +\sum _{j=-\infty}^{\infty}H_{n-j}x_j
\end{equation}
where
$$
H(t)=\sum _{n=-\infty}^{\infty}H_{n}e^{(2\pi int/a)}~.
$$
Since the case in which $H(t)$ and $H_0(t)$ are constant in time is a special case of periodicity, the projection $P$ in (2) cannot be taken as unity here. If in (2), we let $t=na+s$, where $n$ is an integer, we obtain
\begin{equation}
W_{\pm}(s)=\lim_{n\rightarrow \pm \infty}\Theta _{0}^{*n}\Theta^nP
\end{equation}
where $\Theta=U(s+a,s)$ is the operator which takes the system through a complete period, starting at $s$. From this it is clear that one should take
$$
P=P_a(\Theta)~,
$$
the projection onto the absolutely continuous subspace of the unitary operator $\Theta$. Our analysis will connect $P_a(\Theta)$ with the projection
$P_a(\tilde{K})$ which arises  naturally in the ST of $\tilde{K}$.

\medskip

This correspondence between $\tilde{K}$ and $\Theta$ can be illuminated by considering point spectra. For simplicity, take $s=0$ in (5). An eigenvector $\phi(t)$ of $\tilde{K}$ is a solution of
\begin{equation}
\tilde{K}\phi(t)=\left(-i\frac{d}{dt}+H(t)\right)\phi(t)=\lambda \phi(t)
\end{equation}
of period $a$. The functions $\psi(t)=e^{i\lambda t}\phi(t)$ then satisfies $K\psi=0$. Hence, $\psi(t)=U(t,0)\psi(0)$ or
\begin{equation}
U(t,0)\phi(t)=e^{-i\lambda t}\phi(t)~.
\end{equation}
Set $t=a$ and use $\phi(a)=\phi(0)$ to obtain
\begin{equation}
\Theta\phi(0)=e^{-i\lambda a}\phi(0)~.
\end{equation}
This abstract theory is developed in Sections 1 through 5. By way of application, we first obtain (\S 6) an improvement of a result of Schmidt [17] for $H(t)=H_0+V(t)$ with $V(t)$ periodic and of trace class. We then consider Schroedinger's equation (\S 7) and improve the result of Yajima [19] by removing the conditions concerned only with strong solutions of evolutions equations. In both cases, we work in the representation (2).

\medskip

Finally (\S 8), we discuss Schroedinger's equation using the representation (3). In this case, because (3) is closer to the usual representation of stationary ST, we are able to adapt Agmon's method to study the singular spectrum of $U(s+a,s)$ as well. However, although Agmon's theory [20,24] goes through for potentials behaving at infinity like $|x|^{-\alpha}$, $\alpha > 1$, we have required $\alpha > 2$ here. This is because the uniform bound (A.2) of the appendix is required to make the operator-valued multiplication $Q_r(\zeta)$ of \S 8 a bounded operator. A similar problem arises in three-body scattering [23].

\medskip

The most convenient way to extend our results to $\alpha>1$ seems to be to use the `geometrical' method of Enss [22]. We hope to discuss this elsewhere. The methods and results of this section indicate that treating a periodic system is very little different from treating the case of constant Hamiltonians of the same general type.

\bigskip

It should be pointed out that Yajima [19] also employs the representation (3) in the case he considers, although he does not develop the abstract theory as fully as we do. In particular, the correspondence between the absolutely continuous parts of $\tilde{K}$ and $U(s+a,s)$ is noted there.

\medskip

Besides the papers of Schmidt [17] and Yajima [19], there is also that of Davies [2], who has a small perturbation theorem for periodic Hamiltonians. Some other papers on time-dependent Hamiltonians are [1,5,6,7,8,10,11,15,16].

\medskip

One other point deserves mention. In \S 1, we give a new proof of Theorem 1 of [9] by reducing it to von Neumann's theorem on uniqueness of the Schroedinger representation of the {\em CCR}. A proof of von Neumann's theorem can also be based on Theorem 1 plus a little multiplicity theory.

\bigskip

{\bf Notation}. $\qquad$ We denote the integers by $Z$, the reals by $\mathbf R$ and Euclidean $n$-space by ${\mathbf R}^n$. For a normal operator $H$, we denote the absolutely continuous subspace by ${\cal H}_a(H)$, the projection on this space by $P_a(H)$, and the spectrum of $H$ on ${\cal H}_a(H)$ by $\sigma_a(H)$. We define similarly $P_s(H)$ and $P_d(H)$, {\em etc}., corresponding to the singular and point spectra of $H$. We write
$\ell_p({\cal H})$ for the space of ${\cal H}$-valued sequences with $p^{th}$ power summable norms, and $L^p({\mathbf R}^n, {\cal H})$ or $L_p({\mathbf R}^n, {\cal H})$ for the ${\cal H}$-valued $L_p$-space. By $L_{s}^{2}({\mathbf R})$ is meant the weighted space with norm
$$
\int_{{\mathbf R}^n}|f(x)|^2(1+|x|^2)^sdx~.
$$
For a Borel set $S$, $\chi_S$ denotes the characteristic function of $S$.

\bigskip

\S 1. {\bf Evolution groups}. $\quad$ Let ${\cal H}$ be a a separable Hilbert space, and ${\cal H}=L_2({\mathbf R};{\cal H})$. If $A(t)$ is a measurable operator-valued function, then the {\em multiplication} operator ${\cal A}=M(A(t))$ defined by
$$
({\cal A}f)(t)=A(t)f(t)
$$
has norm $||{\cal A}||=ess.sup||A(t)||$. A special case is the {\em scalar multiplication} for $\phi\in L_{\infty}({\mathbf R})$ , with norm
$||N(\phi)||=||\phi||_{\infty}$.
Let $T_{\sigma}$ be the {\em time translation} semigroup
$$
(T_\sigma f)(t)=f(t-\sigma)~.
$$
A unitary group $e^{-i\sigma K}$ on ${\cal H}$ for which $e^{-i\sigma K}T^*_{\sigma}$ commutes with every scalar multiplication $M(\phi)$ is called a
{\em unitary evolution group}. The generator $K$ will be called a {\em (self-adjoint) evolution generator}.

\medskip

{\bf Theorem 1.} $\quad$ {\em A strongly continuous unitary group $e^{-i\sigma K}$ is a unitary evolution group if and only if there is a measurable unitary operator-valued function $U(t)$ such that
\begin{equation}\label{theor1}
e^{-i\sigma K}={\cal U}T_{\sigma}{\cal U}^*
\end{equation}
where ${\cal U}=M(U(t))$}.

\bigskip

A non-unitary version of Theorem 1 was proved in [9, Theorem 1]. In its present form, it can be reduced to von Neumann's theorem on uniqueness of the Schroedinger representation as follows:

\medskip

{\em Proof of Theorem 1}. $\quad$ Clearly (\ref{theor1}) defines a unitary evolution group. Conversely, if $e^{-i\sigma K}$ is such a group, then
$$
e^{-i\sigma K}M(\phi)e^{i\sigma K}=T_\sigma M(\phi)T^*_{\sigma}=M(\phi_\sigma)
$$
where $\phi_\sigma(t)=\phi(t-\sigma)$. Letting $\phi(t)=\chi_S(t)$ for a Borel set $S$ gives
\begin{equation}\label{1.2}
e^{-i\sigma K}E[S]e^{i\sigma K}=E[S+\sigma]
\end{equation}
where $E$ is the spectral measure of the multiplication operator $Qf(t)=tf(t)$. But (\ref{1.2}) is Weyl's form of the canonical commutation relation, so by von Neumann's theorem [3, p. 233], there is a unitary operator ${\cal U}$ from ${\cal H}$ to $L_2({\mathbb R};{\cal H}')$ such that
\begin{equation}\label{1.3a}
{\cal U}^*e^{-i\sigma K}{\cal U}=T'_{\sigma}
\end{equation}
and
\begin{equation}\label{1.3b}
{\cal U}^*E[S]{\cal U}=E'[S]~.
\end{equation}
By (\ref{1.3b}), dim${\cal H}$ = dim${\cal H}'$. Taking ${\cal H}={\cal H}'$ gives $E'[S]=E[S]$ and $T'_\sigma=T_\sigma$.
Equation (\ref{1.3a}) now becomes (\ref{theor1}), and (\ref{1.3b}) says that ${\cal U}$ is a multiplication. \hfill $\blacksquare$

\medskip

The prototype of a unitary evolution group is
\begin{equation}\label{1.4}
e^{-iK\sigma}f(t)=U(t,t-\sigma)f(t-\sigma)
\end{equation}
where $U(t,s)$ is a unitary propagator; that is, $U(t,s)$ is strongly continuous in $(t,s)$, unitary for each $t$ and $s$, and satisfies
$$
U^*(t,s)=U(s,t)
$$
and the Chapman-Kolmogoroff equation
$$
U(t,r)U(r,s)=U(t,s)~.
$$
In this case, one can take $U(t)=U(t,a)$ for any $a$, so that
\begin{equation}\label{1.5}
U(t,s)=U(t)U^*(s)~.
\end{equation}
In the general case, (\ref{1.5}) defines a unique a.e. {\em measurable propagator} associated with the group (\ref{theor1}) [9, p. 318]. If
$U(t,s)$ is the propagator for the evolution equation
$$
i\frac{du}{dt}=H(t)u
$$
on ${\cal H}$, then formally
$$
K=-i\frac{d}{dt}+H(t)
$$
on ${\cal H}$.

\bigskip

\S 2. {\bf Periodic evolution groups}.

\medskip

\noindent The unitary evolution group $e^{-i\sigma K}$ is {\em periodic} with {\em period} $a$ iff $K$ commutes with $T_a$. This corresponds formally to $H(t+a)=H(t)$. Following the standard method for dealing with {\em spatially} periodic Hamiltonians (see, e.g. [18]), we shall transform $K$ to a representation in which $T_a$ is diagonal. We shall for the remainder of this paper take the period $a=1$.

\medskip

For $f\in {\cal H}$, define
$$
({\cal T}f)(n,\theta)=\hat{f}(2\pi n+\theta)
$$
where $n\in {\mathbb Z}$ and $0\leq \theta \leq 2\pi$, and where
$$
\hat{f}(k)=(2\pi)^{-\frac{1}{2}}\int_{-\infty}^{\infty}e^{-ikt}f(t)dt
$$
is the Fourier transform. The map ${\cal T}$ is unitary from ${\cal H}$ onto $L_2((0,2\pi); l_2({\cal H}))$.
In fact, this is just another way of writing $L_2(\hat{{\mathbb R}};{\cal H})$ where $\hat{{\mathbb R}}={\mathbb R}$ is the dual group of ${\mathbb R}$, and we will sometimes write $k=2\pi n+\theta$, in which case ${\cal T}$ is just the Fourier transform. We have, therefore
$$
{\cal T}T_\sigma {\cal T}^*f(k)=e^{-ik\sigma}\hat{f}(k)=e^{-2\pi i n\sigma}e^{-i\theta \sigma}\hat{f}(2\pi n+\theta)
$$
so that
\begin{equation}\label{2.1}
{\cal T}T_\sigma {\cal T}^*=M(e^{- i J\sigma}e^{-i\theta \sigma})
\end{equation}
on $\hat{{\cal H}}$, where $J$ acts on $\ell_2({\cal H})$ according to
\begin{equation}\label{2.2}
(Jx)_n=2\pi n x_n~.
\end{equation}
For $\sigma=1$, we get
\begin{equation}\label{2.3}
{\cal T}T_1 {\cal T}^*=M(e^{-i\theta})
\end{equation}
so that ${\cal T}$ diagonalizes the unitary operator $T_1$.

\bigskip

{\bf Theorem 2}. {\em Let $K$ generate a unitary ev. group with period 1. There exists a unique self-adjoint operator $\tilde{K}$ on $\ell_2({\cal H})$ such that
$$
{\cal T}K {\cal T}^*=M(\tilde{K}+\theta)~.
$$
 The operator $\tilde{K}$ satisfies
\begin{equation}\label{2.4}
Se^{-i\sigma \tilde{K}}=e^{2\pi i \sigma}e^{-i\sigma \tilde{K}}S
\end{equation}
 where $S$ is the bilateral shift on $\ell_2({\cal H})$.

Conversely, given a self-adjoint operator $\tilde{K}$ on $\ell_2({\cal H)}$ satisfying (\ref{2.4}) the operator
$$
K={\cal T}^*M(\tilde{K}+\theta){\cal T}
$$
generates a unitary evolution group of period 1.}\\

Formally, the operator
\begin{equation}\label{2.5}
K=-i\frac{d}{dt}+H_0+V(t)
\end{equation}
on ${\cal H}$ corresponds to
\begin{equation}\label{2.6}
(\tilde{K}x)_n=2\pi n x_n+H_0x_n+\sum_{m=-\infty}^{\infty}\hat{V}_{n-m}x_m
\end{equation}
on $\ell_2(\cal H)$, where
\begin{equation}\label{2.7}
\hat{V}_{n}=(2\pi)^{-1}\int _{0}^{1}e^{-2\pi int}V(t)dt~.
\end{equation}
We shall preface the proof with a few remarks. The interval $0\leq \theta _0<2\pi$ may be regarded as the unit circle ${\mathbb T}_1$. Let
$$
R(\theta_0)f(\theta)=f(\theta-\theta_0)
$$
be the rotation on the circle by the angle $\theta_0$, $\theta_0\in [0,2\pi)$. In terms of its action on $\hat{{\cal H}}$, $R(\theta_0)$ has the somewhat complicated form
\begin{equation}\label{2.8}
R(\theta_0)\hat{f}(2\pi n +\theta)=
\left\{ \begin{array}{ll}
\hat{f}(2\pi n +\theta-\theta_0) & \quad 0\leq \theta < 2\pi -\theta_0\\
\hat{f}(2\pi n +2\pi +\theta-\theta_0) & \quad 2\pi -\theta_0 \leq \theta < 2\pi~.
\end{array} \right.
\end{equation}
If we define the set
$$
B(\theta_0)=\{k: 2\pi n -\theta_0 \leq k < 2\pi n,\, \textrm{for some $n\in {\bf Z}$}\}
$$
and let $B^{c}(\theta_0)$ be its complement, then (\ref{2.8}) may be rewritten on $L_2(\hat{\mathbf{R}},{\cal H})$ as
$$
R(\theta_0)\hat{f}(k)=M(\chi_{B^{c}(\theta_0)}(k)\hat{f}(k-\theta_0)+M(\chi_{B^{c}(\theta_0)}(k)\hat{f}(k-\theta_0 +2\pi)
$$
where $\chi_B$ is the characteristic function of $B$. Equivalently,
\begin{equation}\label{2.9}
R(\theta_0)=M(\chi_{B^{c}(\theta_0)})\hat{T}_{\theta_0}+M(\chi_{B^{c}(\theta_0)})\hat{T}_{(\theta_0-2\pi)}
\end{equation}
where
$$
\hat{T}_\sigma\hat{f}(k)=\hat{f}(k-\sigma)
$$
is translation on ``momentum'' space. Note that
$$
\hat{T}_\sigma={\cal T}M(e^{it\sigma}){\cal T}^{*}~.
$$
We can also express $\hat{T}_{\theta_0}$ in terms of $R(\theta_0)$ and $\hat{T}_{2\pi}$. Multiply (\ref{2.9}) by $\chi_{B^{c}}$ and $\chi_B$ to obtain
$$
M(\chi_{B^{c}})\hat{T}_{\theta_0}=M(\chi_{B^{c}})R(\theta_0) \quad {\rm and} \quad M(\chi_{B})\hat{T}_{\theta_0}\hat{T}_{2\pi}^{*}=
M(\chi_{B})R(\theta_0)~.
$$
Since obviously $\hat{T}_{\theta_0}=[M(\chi_{B^{c}})+M(\chi_{B})]\hat{T}_{\theta_0}$, we obtain
\begin{equation}\label{2.10}
\hat{T}_{\theta_0}=M(\chi_{B^{c}(\theta_0})R(\theta_0)+M(\chi_{B(\theta_0})R(\theta_0)\hat{T}_{2\pi}~.
\end{equation}

Note finally that $\hat{T}_{2\pi}$ is just multiplication on $L_2((0,2\pi),\ell_2({\cal H}))$ by the constant (in $\theta$)
operator $S$, the bilateral shift on $\ell_2({\cal H})$.

\bigskip

{\bf Lemma 2.1.} $\quad$ {\em Let $A=M(A(\theta))$ be a multiplication operator on $L_2((0,2\pi),\ell_2({\cal H}))$. Then $A$ commutes with $T_0$ for every $a$ if and only if $A$ commutes with $\hat{T}_{2\pi}=M(S)$ and with $R(\theta_0)$ for every $\theta_0$.}

\medskip

{\em Proof}. $\quad$ Immediate from (\ref{2.9}) and (\ref{2.10}).

\medskip

{\bf Lemma 2.2.} $\quad$ {\em Let $A$ be an operator on ${\cal K}=L_2(\mathbf{R};{\cal H})$ and $\hat{A}={\cal T}A{\cal T}^{*}$.

\noindent (a) $A$ is a multiplication on ${\cal K}$ if and only if $\hat{A}$ commutes with $\hat{T}_a$ for every $a$.

\noindent (b) If $A$ is a multiplication and commutes with $T_1$, then $\hat{A}$ commutes with $M(e^{-i\theta})$ and $R(\theta_0)$
for every $\theta_0$ and is therefore a contant multiplication on $L_2((0,2\pi),\ell_2({\cal H}))$.}

\medskip

{\em Proof}. $\quad$ If $A$ commutes with every $M(\phi)$, then $\hat{A}$ must commute with ${\cal T}M(e^{-iat}){\cal T}^{*}=\hat{T}_a$.
The converse follows by strong approximation. If, in addition, $A$ commutes with $T_1$, then $\hat{A}$ must commute with ${\cal T}T_1{\cal T}^{*}=
M(e^{-i\theta})$ and hence with any scalar multiplication. Hence, $\hat{A}$ commutes with any scalar multiplication. Hence, $\hat{A}$ commutes with $R(\theta_0)$ by (\ref{2.9}). \hfill $\blacksquare$

\medskip

Finally, note that the operator $J$ defined by (\ref{2.2}) satisfies
\begin{equation}\label{2.11}
e^{iJ\sigma}Se^{-iJ\sigma}=e^{2\pi i\sigma}S
\end{equation}
and
\begin{equation}\label{2.12}
e^{-iJ}=I~.
\end{equation}

{\em \bf Proof of Theorem 2.} $\quad$ Since $e^{-i\sigma K}T_{\sigma}^{*}$ is a multiplication on ${\cal K}$ which comutes with $T_1$,
its transform is (by Lemma 2.2(a)) a constant multiplication on $L_2((0,2\pi),\ell_2({\cal H}))$:
$$
{\cal T} e^{-i\sigma K}T_{\sigma}^{*}{\cal T}^{*}=M(Q_\sigma)~.
$$

Using (\ref{2.1}), one computes that
\begin{equation}\label{2.13}
{\cal T} e^{-i\sigma K}{\cal T}^{*}=M(Q_\sigma)({\cal T}T_\sigma{\cal T}^{*})=M(Q_\sigma e^{-iJ\sigma}e^{-i\theta\sigma})~.
\end{equation}
Therefore,
$$
M(Q_\sigma) e^{iJ\sigma})={\cal T}e^{-i\sigma K}{\cal T}^{*}M(e^{-i\theta\sigma})
$$
is the product of two unitary groups, which commute by (\ref{2.13}). It follows that $Q_\sigma e^{-iJ\sigma}$ is a unitary group on
$\ell_2({\cal H})$. If $\tilde{K}$ is its generator, then
\begin{equation}\label{2.14}
{\cal T} e^{-iK\sigma}{\cal T}^{*}=M(e^{-i\sigma \tilde{K}})M(e^{-i\theta\sigma})=M(e^{-i(K+\theta)\sigma})~.
\end{equation}
Conversely, let $\tilde{K}$ be self-adjoint on $\ell_2({\cal H})$, and define a unitary group on ${\cal K}$ by
$$
e^{-iK\sigma}{\cal T}^{*}M(e^{-i(K+\theta)\sigma}){\cal T}~.
$$
This clearly commutes with $T_1$.

\bigskip
\bigskip

By Lemma 2.2(b), the operator
$$
e^{-iK\sigma}T^{*}_{\sigma}={\cal T}^{*}M\left(e^{-i\tilde{K}\sigma}e^{iJ\sigma}\right){\cal T}
$$
commutes with every scalar multiplication if and only if $M\left(e^{-i\tilde{K}\sigma}e^{iJ\sigma}\right)$ commutes with $\hat{T}_a$ for every $a$. But this operator clearly commutes with every $R(\theta_0)$, since it is multiplication by a function independent of $\theta$. By (\ref{2.10}) this is equivalent to commuting with $\hat{T}_{2\pi}=M(S)$. Hence, $e^{-i\sigma K}$ is a unitary evolution group if and only if
$$
e^{-i\tilde{K}\sigma}e^{iJ\sigma}S=Se^{-i\tilde{K}\sigma}e^{iJ\sigma}
$$
which by (\ref{2.11}) is equivalent to (\ref{2.4}).  \hfill $\blacksquare$

\medskip

As in [9], the condition (\ref{2.4}) can be expressed in other ways.

\bigskip

{\bf Proposition 2.3} $\quad$ {\em A self-adjoint operator $\tilde{K}$ on $\ell_2({\cal H})$ satisfies (\ref{2.4}) if and only if either of the following holds:

(a) $SD(\tilde{K})\subset D(\tilde{K})$ and
\begin{equation}\label{2.15}
\tilde{K}S-S\tilde{K}=2\pi S~.
\end{equation}

(b) For Im $\zeta\neq 0$, the resolvent $\tilde{R}(\zeta)=(\tilde{K}-\zeta)^{-1}$ satisfies
\begin{equation}\label{2.16}
S\tilde{R}(\zeta)-\tilde{R}(\zeta)S=2\pi \tilde{R}(\zeta)~.
\end{equation}
}

The proof follows those of Theorems 3 and 4 of [9].

\bigskip

\S 3. {\bf The period operator}. $\quad$ By the period operator for $K$, we have in mind the operator $\Theta=U(1,0)$, which gives the evolution of the system over a period. If $U(t)=U(t,0)$, then
$$
M(\Theta)={\cal U}^*e^{-iK}T^{*}_{1}{\cal U}~.
$$
Conversely, we have

\medskip

{\bf Proposition 3.1.} {\em The operator $A={\cal U}^*e^{-iK}T^{*}_{1}{\cal U}$ on ${\cal K}$ is multiplication by a constant unitary operator $\Theta$.}

{\em Proof}. $\,$ Since $K$ is an evolution generator, $A$ commutes with every $M(\phi)$. But $A$ also commutes with every $T_\sigma$; for
$$
AT_\sigma={\cal U}^*e^{-iK}T_{1}^{*}({\cal U}T_\sigma)={\cal U}^*e^{-iK}T_{1}^{*}(e^{-i\sigma K}{\cal U})=
{\cal U}^*e^{-iK}e^{-i\sigma K}T_{1}^{*}{\cal U}
$$
$$
=({\cal U}^*e^{-i\sigma K})e^{-iK}T_{1}^{*}{\cal U}=(T_\sigma{\cal U}^*)e^{-iK}T_{1}^{*}{\cal U}=T_\sigma A~.  \qquad \qquad \qquad \qquad \qquad \qquad \qquad \qquad \qquad \qquad \blacksquare
$$

\medskip
Of course, if $U(t)$ is not chosen as $U(t,0)$, $\Theta$ will not be $U(1,0)$, but this is not important.

As a corollary, we have
\begin{equation}\label{3.1}
e^{-iK}={\cal U}M(\Theta){\cal U}^*T_1~.
\end{equation}

\medskip

{\bf Theorem 3.} $\quad$ {\em The operator $e^{-i\tilde{K}}$ is unitarily equivalent to the diagonal operator $\tilde{\Theta}$ on $\ell_2({\cal H})$ defined by
$$
(\tilde{\Theta} x)_n=\Theta x_n~.
$$
}
Compare [19, equation (3.17)].

\medskip

{\em {\bf Lemma 3.2.}} $\quad$ {\em Let $\Phi f(t)=\Theta^{[t]}f(t)$, where $[t]$ denotes the greatest integer in $t$.
Then ${\cal U}_1={\cal U}\Phi^{*}$ commutes with $T_1$.
}

\medskip

{\em Proof}. $\quad$ The equation ${\cal U}T_{1}^{*}Af=T_{1}^{*}{\cal U}f$ translates to $U(t)\Theta f(t+1)=U(t+1)f(t+1)$, which implies that
$$
U(t+1)=U(t)\Theta \qquad {\rm a.e.}
$$
By definition
$$
U_1(t)=U(t)\Theta^{*[t]}
$$
so that $U_1(t+1)=U(t+1)\Theta^{*[t+1]}=U(t)\Theta\Theta^*\Theta^{*[t]}={\cal U}_1(t)$. Hence, ${\cal U}_1T_1=T_1{\cal U}_1$. \hfill $\blacksquare$

\medskip

{\em \bf Proof of Theorem 3.} $\quad$ On ${\cal K}$, we have
$$
e^{-iK}T^{*}_{1}={\cal U}M(\Theta){\cal U}^{*}={\cal U}_1\Phi M(\Theta)\Phi^{*}{\cal U}^{*}_{1}={\cal U}_1M(\Theta){\cal U}^{*}_{1}
$$
since $\Phi$ and $M(\Theta)$ commute. Transformed by ${\cal T}$, this becomes
\begin{equation}\label{3.2}
M(e^{-i\tilde{K}})=\hat{{\cal U}}_1M(\tilde{\Theta})\hat{{\cal U}}^{*}_{1}
\end{equation}
on $L_2((0,2\pi), \ell_2({\cal H})$, where
But ${\cal U}_1$ is a multiplication commuting with $T_1$, so by Lemma 2.2.(b), $\hat{{\cal U}}_1=M(W)$ for some unitary operator $W$ on $\ell_2({\cal H})$. Hence,
$$
M(e^{-i\tilde{K}})=M(W\tilde{\Theta}\tilde{W}^{*})
$$
which implies $e^{-iK}=W\tilde{\Theta}\tilde{W}^{*}$.\hfill $\blacksquare$

\medskip

Theorem 3 is important because it lets us identify the absolutely continuous subspace of $\tilde{K}$ in terms of $\Theta$.

\medskip

{\bf Proposition 3.3.} $\quad$ {\em Define $P_0={\cal T}^{*}M(P_a(\tilde{K})){\cal T}$. Then $P_0={\cal U}M(P_1(\Theta)){\cal U}^{*}$.}

\medskip

{\em Proof}. $\quad$ By (\ref{3.2}), we have
$$
{\cal T}^{*}M(P_a(\tilde{K})){\cal T}={\cal T}^{*}M(P_a(e^{-i\tilde{K}})){\cal T}
={\cal T}^{*}\hat{{\cal U}}_1M(P_a(\Theta))\hat{{\cal U}}_1^{*}{\cal T}
$$
$$
=({\cal U}\Phi^{*})({\cal T}^{*}M(P_1(\tilde{\Theta})){\cal T})(\Phi {\cal U}^{*})=
{\cal U}\Phi^{*}M(P_1(\Theta))\Phi {\cal U}^{*}={\cal U}M(P_a(\Theta)){\cal U}^{*}~. \qquad  \qquad  \qquad \blacksquare
$$

\bigskip
\bigskip

\S 4. {\bf Scattering theory}. $\quad$ We shall now connect the ST for the ev. eq. corresponding to $K$ with the ST of $\tilde{K}$. For any pair of self-adjoint operators $H$ and $H_0$, and any reducing projection $P_0$ of $H_0$, we define the wave operators
\begin{equation}\label{4.1}
\Omega_{\pm}(H,H_0;P_0)=s-\lim_{\sigma\rightarrow \pm \infty} e^{-i\sigma H}e^{i\sigma H_0}P_0~. 
\end{equation}
Let $K$ and $K_0$ generate unitary evolution groups with period 1. By theorem 2,
\begin{equation}\label{4.2}
e^{i\sigma K}e^{-i\sigma K_0}={\cal T}^*M(e^{i\sigma K}e^{-i\sigma K_0}){\cal T}~.
\end{equation}
From Proposition 3.3 and the calculation in [9; p.324], we obtain
\begin{equation}\label{4.3}
{\cal U}^*e^{i\sigma K}e^{-i\sigma K_0} P_0{\cal U}f(t)={\cal U}^*e^{i\sigma K}e^{-i\sigma K_0}{\cal U}M(P_a(\Theta))f(t)=
U^*(t+\sigma)U_0(t+\sigma)P_a(\Theta)f(t)~.
\end{equation}
These formulas lead to the following theorem.

\bigskip

{\bf Theorem 4.}$\quad$ {\em Assume that $\Omega_{\pm}(\tilde{K},\tilde{K}_0; P_a(\tilde{K}_0))$ exists. Then

(a) ${\cal W}_{\pm}=\Omega_{\pm}(K,K_0;P_a(\tilde{K}_0)$ exists and ${\cal U}^*{\cal W}_{\pm}{\cal U}_0=M(W_\pm)$ where $W_\pm$ is partially isometric with initial set $P_a(\Theta_0){\cal H}$.

\medskip

(b) For every $h>0$,
$$
W_{\pm}=s-\lim_{\sigma \rightarrow \pm \infty}h^{-1}\int_{0}^{h}U^*(t+\sigma)U_0(t+\sigma)P_a(\Theta _0)dt~.
$$

Assume, in addition, that the range of $\Omega_{\pm}(\tilde{K}, \tilde{K}_0;P_a( \tilde{K}_0))$ is the absolutely continuous subspace of $\tilde{{\cal K}}$. Then

\medskip

(c) The range of ${\cal W}_\pm$ is
$$
{\cal L}={\cal U}M(P_a(\Theta)){\cal U}^*{\cal K}
$$
and
$$
W_\pm{\cal H}=P_a(\Theta){\cal H}~.
$$

\medskip

(d) The scattering operator ${\cal S}={\cal W}_+{\cal W}^{*}_{-}$ is unitary on ${\cal L}$. The operator ${\cal U}{\cal S}{\cal U}^*$ is multiplication by $S=W_+W^{*}_{-}$, which is unitary on $P_a(\Theta){\cal H}$.}

\medskip

The proof follows the arguments of [19, \S 2] very closely, using (\ref{4.2}),(\ref{4.3}) and Proposition 3.3 where appropriate. We omit it.

\medskip

{\bf Remark}. It is possible to use Abelian limit wave operators in Theorem 4. The only adjustment to be made is that one must then {\em assume} that ${\cal W}_\pm$ are isometric.

\medskip

{\bf Corollary 4.1}.$\quad$ {\em If, in Theorem 4, the function $F(t)=U^*_{0}(t)U(t)$ is continuous in the operator norm, then
$$
W_\pm=s-\lim_{t\rightarrow \pm \infty}U^{*}_{0}(t)U(t)P_a(\Theta_0)~.
$$
}

\medskip

{\em Proof}. $\quad$ According to [9, p. 324], it suffices to show that $F(t)$ is {\em uniformly} continuous on $(-\infty,\infty)$. Let $0\leq t-s\leq 1$. Then from $U(t+1)=U(t)\Theta$ we get $U(t+n)=U(t)\Theta^n$ and hence
$$
F(t)-F(s)=\Theta^{*[s]}[F(t-[s])-F(s-[s])]\Theta^{[s]}~.
$$
But $0\leq s-[s]\leq t-[s]\leq 2$, and $F(t)$ is uniformly continuous on $0\leq t\leq 2$.

If $U_0(t)$ is {\em strongly} continuous, then in the corollary, $U(t)$ is also strongly continuous. If, as we may, we take $U(0)=I$, then we find that
$$
\Theta=U(1,0) \equiv U(n+1,n)
$$
and that
$$
W_\pm=s-\lim_{t\rightarrow \infty}U(0,t)U_0(t,0)P_a(U_0(1,0))~.
$$
Transforming by $U(s,0)$ gives that
$$
W_\pm=s-\lim_{t\rightarrow \infty}U(s,t)U_0(t,s)P_a(U_0(s+1,s))~.
$$
is isometric from $P_a(U_0(s+1,s)){\cal H}$ onto $P_a(U(s+1,s)){\cal H}$, and that $S(s)=W_+(s)W^{*}_{-}(s)$ is unitary
on $P_a(U_0(s+1,s)){\cal H}$.

\medskip

\S 5. {\bf Representations on periodic functions}. $\quad$ The mapping ${\cal S}$ which takes the sequence $x=\{x_n\}$ in $\ell_2({\cal H})$ into the function
$$
({\cal F}x)(t)=\hat{x}(t)=\sum_{n=-\infty}^{+\infty} e^{2\pi i nt}x_n
$$
maps $\ell_2({\cal H})$ unitarily onto the space $\hat{{\cal H}}$ of locally square-integrable ${\cal H}$-valued functions of period 1. If $U(t,s)$ is strongly continuous it is easy to see that
\begin{equation}\label{5.1}
{\cal F}e^{-i\tilde{K}\sigma}{\cal F}^*f(t)=U(t,t-\sigma)f(t-\sigma)~.
\end{equation}
If $f$ has period 1, then so does the right side, because $U(t+n, s+n)=U(t,s)$. Equation (\ref{5.1}) is formally the same as Equation ({\ref{1.4}) for the group $e^{-iK\sigma}$, except that $f$ is periodic rather than square-integrable. Hence, the operator
$$
\hat{K}={\cal F}\hat{K}{\cal F}^*
$$
is formally the same as $K$, but it acts on periodic functions. Thus, formally,
$$
\hat{K}=-i\frac{d}{d t}+H(t)
$$
with the periodic boundary condition $u(0)=u(1)$. Since we will usually define $\hat{K}$ by means of its resolvent or unitary group, it will not be necessary to make this more precise.

The following analogue of Proposition 2.3 holds.

\medskip

{\em {\bf Proposition 5.1.} $\quad$ Let $\hat{K}$ be self-adjoint on $\hat{{\cal K}}$, and $M=M(e^{2\pi i t})$. The operator $\hat{K}={\cal F}\hat{K}{\cal F}^*$ satisfies (\ref{2.4}) if and only if any one of the following holds:

(a) $MD(\hat{K})\in D(\hat{K})$ and
\begin{equation}\label{5.2}
\hat{K}M-M\hat{K}=2\pi M~.
\end{equation}

(b) For Im $\xi \neq 0$, the resolvent $\hat{R}(\zeta)=(\hat{K}-\zeta)^{-1}$ satisfies
\begin{equation}\label{5.3}
M\hat{R}(\zeta)-\hat{R}(\zeta)M\hat{R}(\zeta)~.
\end{equation}

(c) One has for all real $\sigma$
\begin{equation}\label{5.4}
Me^{-i\sigma \hat{K}}=e^{2\pi i \sigma}e^{- i \sigma \hat{K}}M~.
\end{equation}
}
The proof is immediate, because ${\cal F}S{\cal F}^*=M$.

This proposition makes it possible to construct operators $\hat{K}$ by constructing the corresponding $\tilde{K}$. This is done in \S 6 and \S 7, and by Yajima [9].
A formula for the resolvent of $\hat{K}$ is easily derived from (\ref{5.1}). For the special case $H(t)\equiv H$, it is given by equation (\ref{6.4}) below. Compare [19].

\bigskip
\bigskip

\S 6. {\bf Trace class perturbations}. $\quad$ The following generalizes a result of E.J.P. Georg Schmidt [17]. Let $||V||_1$ and $||V||_2$ denote the trace and Schmidt norms of $V$, respectively.

\bigskip

{\em {\bf Theorem 5}. $\quad$ Let $V(t)$ be a measurable self-adjoint trace-class-valued function on $-\infty < t < \infty$. Assume that $V(t)$ has period 1 and that
\begin{equation}\label{6.1}
\int_{0}^{1}||V||_1dt < \infty~.
\end{equation}
Let $H$ be self-adjoint, and let $U(t,s)$ be the propagator for $H(t)=H+V(t)$.
Then
$$
W_\pm(s)=st-\lim_{t\rightarrow \pm \infty} U^{*}(t,s)e^{-iH(t-s)}P_{ac}(H)
$$
exist, and are partially isometric with ranges equal to ${\cal H}_{ac}(U(s+1,s))$.
}

\medskip

The propagator $U(t,s)$ is obtained in the usual way by iteration of the Integral Equation
\begin{equation}\label{6.2}
U(t,s)=U_0(t,s)-i\int_{s}^{t}U_0(t,y)V(y)U(y,s)dy
\end{equation}
where $U_0(t,s=e^{-iH_0(t-s)}$. There is no problem for bounded operators with locally integrable norm.
The ratio $F(t,s)=U_0(s,t)U(t,s)$ satisfies
\begin{equation}\label{6.3}
F(t,s)=I-i\int_{s}^{t}U_0(s,y)V(y)U(y,s)dy
\end{equation}
and is therefore uniformly continuous in norms for $|t|\leq 2$, $|s|\leq 2$. Thus, it suffices to consider the $(\hat{K}hat{K}_0)$} scattering problem. In the Fourier series representation of \S 5, this means considering
$$
\hat{K}=-i\frac{d}{dt}+H+V(t)
$$
on functions on $0\leq t\leq 1$ with the periodic boundary condition $u(0)=u(1)$.

The operator
$$
\hat{K}_0=-i\frac{d}{dt}+H
$$
has resolvent
\begin{equation}\label{6.4}
\hat{R}_0(\lambda)f(t)=i\int_{0}^{t}e^{-i\lambda(s-t)}e^{iH(s-t)}f(s)ds+i 
\frac{\int_{0}^{1}e^{-i\lambda(s-t)}e^{iH(s-t)}f(s)ds}{e^{-i\lambda}e^{iH}-1}
\end{equation}
which for Im $\lambda < 0$, has the series expansion
$$
\hat{R}_0(\lambda)f(t)=i\int_{0}^{t}e^{-i\lambda(s-t)}e^{iH(s-t)}f(s)ds
-i\sum_{n=1}^{\infty}
\int_{0}^{1}e^{-i\lambda(s-t+n)}e^{iH(s-t+n)}f(s)ds
$$
where the first terms of the series have been combined with the first term of (\ref{6.4}).

Let $A$ and $B$ be multiplications by functions $A(t)$ and $B(t)$ of period 1, and assume that
\begin{equation}\label{6.5}
M=\int_{0}^{1}||A(t)||_{2}^{2}dt < \infty
\end{equation}
and
\begin{equation}\label{6.6}
N=\int_{0}^{1}||B(t)||^{2}dt < \infty~.
\end{equation}
For $\eta=-{\rm Im} \lambda >0$, one has
\begin{equation}\label{6.7}
A\hat{R}_0(\lambda)Bf(t)=i\int_{t}^{1}e^{-i\lambda(s-t)}A(t)e^{iH(s-t)}B(s)f(s)ds
-i\sum_{n=1}^{\infty}\int_{0}^{1}e^{-i\lambda(s-t+n)}A(t)e^{iH(s-t+n)}B(s)f(s)ds~.
\end{equation}
The squared Schmidt norm of the first term is
$$
\int_{0}^{1}\int_{t}^{1}e^{-2\eta\lambda(s-t)}||A(t)e^{iH(s-t)}B(s)||_{2}^{2}dsdt\leq
\int_{0}^{1}\int_{0}^{1}e^{-2\eta\lambda(s-t)}||A(t)||_{2}^{2}B(s)||^{2}dsdt=O(1)
$$
as $\eta\rightarrow \infty$. The squared Schmidt norm of the n$^{th}$ term of the series in (\ref{6.7}) is
$$
\int_{0}^{1}\int_{0}^{1}e^{-2\eta(s-t+n)}||A(t)e^{iH(s-t)}B(s)||_{2}^{2}dsdt
$$
$$
\leq e^{-2\eta(n-1)}\int_{0}^{1}\int_{0}^{1}||A(t)||_{2}^{2}||B(s)||^{2}e^{-2\eta(s-t+1)}dsdt=O(e^{-2\eta(n-1)})
$$
uniformly in $n$ as $\eta\rightarrow \infty$. Hence, $||A\hat{R}_0(\lambda)B||_2$ is finite for $\eta > 0$ and $O(1)$ as $\eta \rightarrow \infty$.

We can now use the factorization method to define $\tilde{K}$ [12,14,p.148]. Both $A(t)=|V(t)|^{\frac{1}{2}}$ and $B(t)=|V(t)|^{\frac{1}{2}}$sgn$V(t)$ satisfy (\ref{6.5}) and (\ref{6.6}). It follows that $Q(\lambda)=A\hat{R}_0(\lambda)B$ is bounded, and that $I+Q(\lambda)$ is invertible for Im $\lambda$ large and hence for all nonreal $\lambda$. Also, $A\hat{R}_0(\lambda)$ and $B\hat{R}_0(\lambda)$ are Hilbert-Schmidt, as can be seen from the above estimates since $B(t)\equiv I$ satisfies (\ref{6.6}). Following [12,14], one obtains a self-adjoint $\hat{K}$ satisfying (\ref{2.16}), with resolvent
\begin{equation}\label{6.8}
\hat{R}(\lambda)=\hat{R}_0(\lambda)-[B\hat{R}_0(\lambda)]^*[I+Q(\lambda)]^{-1}A\hat{R}_0(\lambda)~.
\end{equation}

\medskip

Since $A\hat{R}_0(\lambda)$ and $B\hat{R}_0(\lambda)$ are Hilbert-Schmidt, the difference $\hat{R}(\lambda)-\hat{R}_0(\lambda)$ is trace class. The result now follows from a classical theorem of Birman. See [13,p. 548] for a proof using the invariance principle.

\bigskip

{\bf Remark}. $\quad$ It is possible to base a proof of Theorem 5 on an old result of Birman and Krein [21], which asserts the existence of complete wave operators for any two unitaries which differ by a trace-class operator. From (\ref{6.1}) and (\ref{6.2}), $U(t,s)$ and $U_0(t,s)$ are such operators, so
$$
\lim U(0,t)U_0(t,0)
$$
exists as $t$ tends to infinity through {\em integer} values. One must then extend this to all $t$ ({\em cf}. [17]). This is likely to be simpler than the proof above, but it does not illustrate our method.

\bigskip

\S 7. {\bf Schroedinger's equation}. $\quad$ We shall consider the equation
\begin{equation}\label{7.1}
i\frac{\partial \psi(x,t)}{\partial t}=-\Delta_n\psi(x,t)+q(x,t)\psi(x,t)~,
\end{equation}
on $\mathbb{R}^\nu$, $\nu\geq 3$, where $q(x,t)$ has period 1 in $t$. We shall extend the result of Yajima [19] by eliminating the smoothness assumptions in $t$ which are only used to obtain strong solutions of the corresponding evolution equation. Many estimates are similar to those of [19], so we will be rather sketchy. We shall assume:

\medskip

{\bf Hypothesis A}. $\quad$ {\em Let $q(x,t)$ be real-valued and measurable on $\mathbb{R}^{\nu+1}$, $\nu\geq 3$, and periodic in $t$ with period 1.
Let
$$
V_p(t)=\left(\int_{R^\nu}|q(x,t)|^pdx\right)^{\frac{1}{p}}~.
$$
Assume that there exist $p$ and $q$, $1\leq q < \nu/2<p\leq \infty$, such that $V_p\in L_\sigma[0,1]$ for some $\sigma$, $(2p/2p-\nu)<s\leq \infty$,
and $V_q\in L_r[0,1]$ for some $r$, $2\leq r\leq \infty$.
}

\medskip

As in the preceding section, we shall use the Fourier series representation of $\hat{K}$ on periodic functions on $0\leq t\leq 1$. For $\hat{K}_0$, take
$$
\hat{K}_0=-i\frac{d}{dt}+H
$$
with $H=-\Delta_\nu$ on ${\cal H}=L_2(\mathbb{R}^\nu)$. Let $A$ and $B$ be multiplication by $a(x,t)=|q(x,t)|^{\frac{1}{2}}$ and $b(x,t)=|q(x,t)|^{\frac{1}{2}}$sgn $q(x,t)$, and set $Q(\lambda)=A\hat{R}_0(\lambda)B$.

\medskip

{\em {\bf Lemma} $\quad$ The operators $A$ and $B$ are $\hat{K}_0$-smooth, and $Q(\lambda)$ is compact, continuous up to the axis in norm, and tends to zero in norm as $|{rm Im} \lambda|$ tends to infinity.}

\medskip

{\em Proof}. $\quad$ According to Kato [12; Equation (6.9)], the norm of $A(t)e^{iH(s-t)}B(s)$ does not exceed $(4\pi |s-t|^{-\nu/2p}[V_p(t)V_p(s)]^{1/2}$. In order to estimate $Q(\lambda)$, combine the first term of the series in (\ref{6.2}) with the indefinite integral, and use Kato's estimate to obtain:
$$
|A\hat{R}_0(\lambda)Bf(t)|\leq C\int_{t}^{2}[V_p(t)V_p(s)]^{1/2}(s-t)^{-\nu/2p}|f(s)|ds+C\sum_{k\geq 2}\int_{0}^{1}[V_q(t)V_q(s)]^{1/2}(s-t+k)^{-\nu/2q}|f(s)|ds
$$
The operators appearing here are compositions of multiplications and convolutions, and can be estimated by Holder's and Young's inequalities, exactly as in [9; Lemma, p. 333]. Then the $L_2$-norm of the first term does not exceed
$$
C||V_p||_{\sigma}\left(\int_{0}^{2}s^{-\nu \sigma'/2p}ds\right)^{1/\sigma'}||f||_2
$$
while that of the second term does not exceed
$$
C||V_q||_{r}||f||_2\sum_{k\geq 2}\left(\int_{0}^{2}(s+k-1)^{-\nu\tau'/2q}ds\right)^{1/\tau'}\leq
2^{1/\tau'}C||V_q||_{r}||f||_2\sum_{k\geq 2}(k-1)^{-\nu/2q}< \infty~.
$$
This gives a norm estimate of $Q(\lambda)$ uniformly in $\lambda$. The proof of continuity up to the axis, and the vanishing of the norm at infinity requires only to retain the appropriate exponential factors $e^{-i\lambda(s-t)}$in the estimates ({\em cf}. [9, p. 331]). Compactness may be proved by Hilbert-Schmidt approximation, as in [9].

\medskip

It follows that (\ref{6.6}) (reinterpreted) defines a self-adjoint $\hat{K}$ whose resolvent is easily seen to satisfy (\ref{5.3}).
Let $U(t,s)=U(t)U^*(s)$ be the (measurable) propagator for the corresponding operator $K$ on ${\cal H}$, and $U_0(t,s)=e^{-iH(t-s)}$.

\bigskip

{\em {\bf Theorem 6}. $\quad$ The propagator $U(t,s)$ is strongly continuous. The limits
$$
W_{\pm}(s)=s-\lim_{t\rightarrow \pm \infty} U_0(s,t)U(t,s)P_a((U(s+1,s)) 
$$
exist, and
$$
W_{\pm}(s)^*=s-\lim_{t\rightarrow \pm \infty} U(t,s)U_0(s,t)~.
$$
}

\medskip

{\em Proof}. $\quad$ The scattering theory for $\hat{K}$ and $\hat{K}_0$ may be obtained by inverting $I+Q(\lambda)$ on the axis, using the Kato-Kuroda lemma [14, Lemma 4.20]. Thus, it suffices to prove strong continuity of $U(t,s)$ and uniform continuity of the quotient $F(t)$ of Corollary 4.1. This will follow from the fact that $q(x,t)$ satisfies {\em locally} in $t$ the condition of [9, \S 7], so that if we cut off our periodic $q(x,t)$ for large $t$, we get an evolution equation satisfying that condition. The propagator of the new equation is strongly continuous by the result of [9]. But for $t,s$ in some interval, $U(t,s)$ depends only on the values of $q(x,t)$ for $t$ in that interval, so the propagators should agree locally.

\medskip

We shall make this argument precise. If $K$ and $K_n$ are self-adjoint, we say that $K_n\rightarrow K$ {\em strongly in the generalized sense}
(sgs) if and only if
\begin{equation}\label{7.2}
(K_n-\zeta )^{-1}\rightarrow (K-\zeta)^{-1}
\end{equation}
strongly for Im $\zeta \neq 0$. [13, p. 427].

\medskip

{\em {\bf Lemma}. $\quad$ If $K_n$ generates a unitary evolution group for each $n$, and $K_n\rightarrow K$ (sgs), then $K$ also generates a unitary evolution group. If $U_n(t,s)$ and $U(t,s)$ are the corresponding propagators, then
\begin{equation}\label{7.3}
\lim \int _S \langle U_n(t,s)x,y\rangle dt ds=\int _S \langle U(t,s)x,y\rangle dt ds
\end{equation}
for every $x,y \in {\cal H}$ and every Borel subset S of $\mathbb{R}^2$.}

\medskip

{\em Proof}. $\quad$ To see that $K$ generates a unitary evolution group, use strong convergence and [9, Theorem 3, p. 322]. For Im $\zeta > 0$, (\ref{7.2}) means in the present case that
$$
\phi_n(t)=ie^{i\zeta t}\int_{-\infty}^{t}D_n(t,s)e^{-i\zeta s}f(s)ds
$$
tends to zero in the $L_2$-norm, where $D_n(t,s)=U_n(t,s)-U(t,s)$. Let $f(s)e^{i\zeta s}\xi_E(s)x$ and $\psi(t)=e^{i\zeta s}\xi_F(s)y$,
where $\xi_E$ and $\xi_F$ are the characteristic functions of the bounded Borel sets $E$ and $F$. Then
$$
0=\lim \langle \phi_n,\psi\rangle=\lim \int\!\!\!\int_{Q\bigcap G}\langle D_n(t,s)x,y\rangle ds dt
$$
where $Q=E\times F$ and $G=\{(t,s)|t>s\}$. Considering Im $\zeta<0$ gives the same result with $G$ replaced by its complement $G^c$, so that we obtain:
$$
\lim \int\!\!\!\int_{Q}\langle D_n(t,s)x,y\rangle ds dt
$$
for all finite rectangles $Q$, which is equivalent to (\ref{7.3}). This proves the lemma. \hfill $\blacksquare$

\medskip

To prove Theorem 6, approximate $q(x,t)$ in the appropriate norms by a sequence $q_n(x,t)$ of bounded functions of period 1 in $t$. Let $\xi(t)\in C_{c}^{\infty}(\mathbb{R}^n)$ be identically 1 for $-2\leq t \leq 2$, and $0\leq \xi(t)\leq 1$. Define
$$
K_n=K_0+q_n
$$
and
$$
L_n=K_0+\xi q_n~.
$$
Then $K_n\rightarrow K$ (sgs) and $L_n\rightarrow K_0+\xi q$ (sgs). The operators $K_n$ and $L_n$ are {\em bounded} perturbations of
$$
\hat{K}_0=-i\frac{d}{dt}+H
$$
by multiplications, and therefore generate unitary evolution groups. The corresponding propagators $U_n(t,s)$ and $W_n(t,s)$ are strongly continuous, and are obtained by iterating integral equations like (\ref{6.2}). Hence, clearly $U_n(t,s)=W_n(t,s)$ for $-2\leq t, \, s\leq 2$.
Hence, by the lemma, if $S\subset [-2,2]\times [-2,2]$, then
$$
\int_S\langle U(t,s)x,y\rangle dt ds=\lim \int_S\langle U_n(t,s)x,y\rangle dt ds
$$
$$
=\lim \int_S\langle W_n(t,s)x,y\rangle dt ds=\lim \int_S\langle W(t,s)x,y\rangle dt ds~.
$$
But $L$ satisfies the condition of [9,\S 4], so $W(t,s)$ is strongly continuous. Hence, after redefinition on a set of measure zero, $U(t,s)=W(t,s)$ for $-2\leq t$, $s\leq 2$. But $U$ has period 1, and so $U(t,s)$ is strongly continuous. Moreover, [9, p. 329]
$$
F(t)=W^*(t)W_0(t)-U^*(t)U_0(t)
$$
is continuous in norm for $-2\leq t\leq 2$. The theorem now follows. \hfill $\blacksquare$

\bigskip

\S 8. {\bf Schroedinger's equation, II}. $\quad$ We shall now consider equation (\ref{7.1}) with
$$
q(x,t)=\sum _{n=-\infty}^{+\infty} e^{2\pi i nt}V_n(x)~.
$$
For simplicity, we shall consider only dimension $\nu=3$. We shall assume that for some $\alpha>2$,
$$
V_n(x)=(1+|x|^2)^{-\alpha/2}W_n(x)
$$
where $W_n\in L^p(\mathbb{R}^3)$ for some $p> 3/2$, and that
\begin{equation}\label{8.1}
B=\sum_n \beta_n<\infty
\end{equation}
where $\beta_n=\beta(V_n)$ denotes the $L^p+L^\infty$ norm of $W_n$. We also assume that $q(x,t)$ is real, which corresponds to
\begin{equation}\label{8.2}
\overline{V_n(x)}=V_{-n}(x)~.
\end{equation}
The series
$$
W(x,t)=\sum_n e^{2\pi i nt}W_n(x)
$$
converges absolutely in $(L^p+L^{\infty})$-norm. Hence
$$
q(x,t)=(1+|x|^2)^{-\alpha/2}W(x,t)
$$
where $W(x,t)$ is uniformly bounded (in $t$) in the $L^p+L^\infty$ norm (in $x$). It follows easily that $q(x,t)$ satisfies the conditions of
\S 7. We want to show that by working with $\tilde{K}$ in $\ell_2({\cal H})$, one can go further and describe the singular part of $\Theta$. Let
$$
(\tilde{K}_0x)_n=(H+2\pi n)x_n
$$
and
$$
\tilde{K}=\tilde{K}_0+V
$$
where $H=-\Delta$ on ${\cal H}=L^2(\mathbb{R}^3)$ and
$$
(Vx)_n=\sum_kV_{n-k}x_k~.
$$

\medskip

{\em {\bf Proposition 8.1}. $\quad$ One has $D(V) \subset D(\tilde{K}_0)$. The operator $\tilde{K}$ is self-adjoint with $D(\tilde{K})=D(\tilde{K}_0)$ and satisfies (\ref{2.15}).}

\medskip

{\em Proof}. $\quad$ The operator $Q(\zeta)=V(\tilde{K}_o-\zeta)^{-1}$ is given by
\begin{equation}\label{8.3}
(Q(\zeta)x)_n=\sum_kV_{n-k}R(z-2\pi k)x_k
\end{equation}
where $R(\zeta)=(H-\zeta)^{-1}$ on ${\cal H}=L^2(\mathbb{R}^3)$. From (\ref{A.1}), one obtains
$$
|Q(\zeta)x_n|\leq C\left(\sum_k \beta_{n-k}|x_k| \right)O(\eta^{-1/4})
$$
where $\eta={\rm Im}\,\zeta$. Using that $\ell_1*l_2\subseteq l_2$, we have
$$
|Q(\zeta)x_n|\leq C B|x|O(\eta^{1/4})~.
$$
Hence, $Q(\zeta)$ tends to zero as $\eta$ tends to infinity, and, $V$ is $\tilde{K}_0$ $\epsilon$-bounded. Symmetry of $\tilde{K}_0$
follows from (\ref{8.2}). Equation (\ref{2.15}) can be verified directly.

\medskip

The results of \S 7 give:

\bigskip

{\em {\bf Theorem 6'}. $\,$ The operators $\Omega_{\pm}(\tilde{K},\tilde{K}_0)$ exist, and their range is the absolutely continuous subspace of $\tilde{K}$.}

\medskip

This result can be obtained by working directly with $\tilde{K}$ and $\tilde{K}_0$. We sketch the proof, since the estimates involved will be used below. One considers $Q(\zeta)$ as an operator on the space ${\cal H}=l_2(L_{s}^{2}(\mathbb{R}^3))$ for an appropriate $s>1$. Using Formula (\ref{8.3}), and the estimates (\ref{A.2}) and (\ref{A.3}), one shows that $Q(\zeta)$ is bounded on ${\cal H}$ and continuous up to the axis in norm.

For compactness, observe that for fixed $r$, the operator
$$
(Q_r(\zeta)x)_n=a_r(\zeta +2\pi r-2\pi n)x_n
$$
is compact on ${\cal H}$ since $a_r(\zeta)$ is compact and (\ref{A.4}) holds. But
$$
Q(\zeta)=\sum_r Q_r(\zeta)
$$
and the sum is norm convergent since, by (\ref{A.2}), $||Q_r(\zeta)||_{{\cal H}}\leq C\beta_r$. By a lemma of Kato and Kuroda [14, Lemma 4.20],
$I+Q(\lambda\pm i0)$ is invertible for a.e. $\lambda$. A spectral form for $\tilde{K}_0$ is
\begin{equation}\label{8.4}
F(\lambda;x)=\sum_kf_0(\lambda-2\pi k; x_k)
\end{equation}
where$x\in l_2({\cal H})$ and $f_0$ is the spectral form of (\ref{A.5}). The result now follows from the stationary theory of [14].

\medskip

We shall now study the singular spectrum of $\tilde{K}$.

\medskip

{\em {\bf Theorem 7}. $\quad$ The singular spectrum $\sigma_s(\tilde{K})$ of $\tilde{K}$ consists of eigenvalues of finite multiplicity, plus perhaps the set $2\pi Z$ of thresholds. The only possible accumulation points of $\sigma_s(\tilde{K})$ are the thresholds. If $\lambda\in \sigma_s(\tilde{K})$, then $\lambda+2\pi n \in \sigma_s(\tilde{K})$ for every $n\in Z$.}

\medskip

{\em Proof}. $\,$ Let $\lambda \notin 2\pi Z$ be a point where $I+Q(\lambda+i0)$ is not invertible, and let $\phi \in {\cal H}$ satisfy $\phi+
Q(\lambda+i0)\phi=0$, that is,
 \begin{equation}\label{8.5}
 \phi_n+\sum_kV_{n-k}R(\lambda+i0-2\pi k)\phi_k=0
 \end{equation}
 for all $n\in Z$. Define a sequence $\psi=(\psi_n)$ with the components
  \begin{equation}\label{8.6}
  \psi_n=R(\lambda+i0-2\pi n)\phi_n~.
   \end{equation}
 We claim that $\psi \in l_2({\cal H})$. For, as usual, [3, Proof of Proposition (3.5)]
 $$
 F(\lambda,\phi)=0
 $$
 which by (\ref{7.4}) implies that
 $$
 f_0(\lambda-2\pi k, \phi_k)=0
 $$
 for every $k$. According to (\ref{A.6})
  \begin{equation}\label{8.7}
  ||\psi_k||\leq C(\delta)||\phi_k||_s
  \end{equation}
  where $C(\delta)$ is independent of $K$ and $\lambda$, as long as dist$(\lambda, 2\pi Z)\geq \delta>0$. It follows that $\psi\in l_2({\cal H})$ and
 \begin{equation}\label{8.8}
  ||\psi||\leq C(\delta)||\phi||_{{\cal H}}~.
  \end{equation}
  By (\ref{8.5}),$\phi=-V\psi$, so $\phi\neq 0$ implies $\psi \neq 0$.
  By definition of $\psi$, $\psi \in D(\tilde{K}_0)\subset D(V)$ and
  \begin{equation}\label{8.9}
  (\tilde{K}_0-\lambda)\psi=\phi =-V\psi
  \end{equation}
  by (\ref{8.5}) and (\ref{8.6}). Hence, $\lambda$ is an eigenvalue of $\tilde{K}_0$.

  Conversely, suppose that $\psi \in D(\tilde{K})$ and $\tilde{K}\psi=\lambda \psi$. Define $\phi=-V\psi(\tilde{K}_0-\lambda)\psi$. An easy calculation gives
  $$
  \phi+Q(\lambda+i\epsilon)\phi=i\epsilon Q(\lambda+i\epsilon)\psi
  $$
  for $\epsilon >0$. From (\ref{A.7}), it follows that $Q(\lambda+i\epsilon)$ maps $\ell_2({\cal H})$ into ${\cal H}=l_2(L_{s}^{2})$ for $\epsilon>0$  . Hence, $\phi\in {\cal H}$ and $\phi+Q(\lambda+i0)\phi$ exists strongly in ${\cal H}$. Let $\chi$ be a finite sequence with
  $\chi_n\in C_{c}^{\infty}$. Such $\chi$'s are dense in ${\cal H}$. We have
$$
i\epsilon \langle Q(\lambda+i\epsilon)\psi,\chi\rangle =\langle i\epsilon(\tilde{K}_0-\lambda-i\epsilon)^{-1}\psi , V\chi\rangle~.
$$
But $V\chi\in l_2({\cal H})$ and
\begin{equation}\label{8.10}
\lim_{\epsilon \rightarrow 0}i\epsilon(\tilde{K}_0-\lambda-i\epsilon)^{-1}\psi=0
\end{equation}
weakly in $\ell_2({\cal H})$ since ker $\tilde{K}_0=(0)$. (This is true for any self-adjoint operator.) It follows that the right side of (\ref{8.9}) tends to zero weakly and hence that
$$
\phi+Q(\lambda+i0)\phi=0~.
$$
Hence, the correspondence between $\phi\in {\rm ker}[I+Q(\lambda+i0)]$ and $\psi \in {\rm ker}(\tilde{K}-\lambda)$ is one-to-one, so that $\lambda$ has finite multiplicity equal to
$$
{\rm dim}\,{\rm ker}\, (I+Q(\lambda+i0))~.
$$
Finally, we claim that there are no accumulation points of $\sigma_s(\tilde{K})$ except possibly $2\pi Z$. Suppose that $\lambda_r$ are distinct points of $\sigma_s(\tilde{K})$, $\lambda_r\rightarrow \lambda_0\notin 2\pi Z$. We can then assume that dist$(\lambda_r, 2\pi Z)\geq \delta>0$. Let
 \begin{equation}\label{8.11}
 \phi^{(r)}=Q(\lambda_r+i0)\phi^{(r)}
 \end{equation}
 with $||\phi^{(r)}||_s=1$. The right side of (\ref{8.11}) ranges over a precompact set, so by passing to a subsequence, we may assume that $\phi^{(r)}\rightarrow \phi^{(0)}$ in ${\cal H}$-norm. Then
 $$
 \phi^{(0)}=-Q(\lambda_0+i0)\phi^{(0)}
 $$
 and $||\phi^{(0)}||_s=1$. If $\psi^{(r)}$ are the corresponding eigenfunctions of $\tilde{K}$, then $\psi^{(0)}\neq 0$ since $\phi^{(0)}\neq 0$, and
  \begin{equation}\label{8.12}
  ||\psi^{(r)}||\leq C(\delta)||\phi^{(r)}||_s=C(\delta)~.
 \end{equation}

 \bigskip

Let $y=(y_n)$ be a finite sequence such that for every $n$, $y_n\in C_{c}^{\infty}(\mathbb{R}^3)$ and $\hat{y}_n(k)=0$ on the set
$$
\{k: |k^2-\lambda_0+2\pi n|<\epsilon\}
$$
for some positive $\epsilon$, depending on $y$. Such $y'$s are dense in $\ell_2({\cal H})$ (although not in ${\cal H}$). For large $r$,
$|\lambda_r-\lambda_0|<\epsilon$. Hence, $R(\lambda_r-i0-2\pi n)y_n$ is in ${\cal H}$ and converges strongly in ${\cal H}$ to $R(\lambda_0-i0-2\pi n)y_n$ as $r$ tends to infinity. Therefore,
$$
u_r=(\tilde{K}_0-\lambda_r+i0)^{-1}y
$$
converges strongly to
$$
u_0=(\tilde{K}_0-\lambda_0+i0)^{-1}y
$$
and we have
$$
\langle \psi^{(0)},y\rangle=\langle \phi^{(0)},u_0\rangle=\lim \langle \phi^{(r)},u_r\rangle=\lim \langle \psi^{(r)},y\rangle
$$
Since, by (\ref{8.12}), the norms of $\psi^{(r)}$ are uniformly bounded, $\psi^{(r)}\rightarrow \psi^{(0)}$ weakly in $\ell_2({\cal H})$.
But $\psi^{(r)}$ is orthogonal to $\psi^{(0)}$, since they are eigenvectors of $\tilde{K}$ for different eigenvalues. This implies $\psi^{(0)}=0$, which contradicts $\psi^{(0)}\neq 0$.

That $\sigma_s(\tilde{K})+2\pi Z=\sigma_s(\tilde{K})$ now follows immediately from (\ref{2.15}). \hfill $\blacksquare$

\medskip

{\em {\bf Corollary}. 
${\cal H}_s(\Theta)$ is spanned by the eigenvectors of $\Theta$. The eigenvalues of $\Theta$ can accumulate only at 1.}

\medskip

\noindent In general, as indicated in the Introduction, one cannot expect the eigenvalues of $\Theta$ to have finite multiplicity.

\bigskip
\bigskip

{\bf Appendix}. $\quad$ We collect here some well-known estimates for the two-body problem in dimension $\nu=3$.

\medskip

\noindent Let $H=-\Delta$ on $L^2(\mathbb{R}^3)$ and $R(\zeta)=(H-\zeta)^{-1}$. Let $V(x)$ be measurable, and perhaps complex-valued, on $\mathbb{R}^3$.
Fix a number $\alpha >2$, and define $\beta(V)$ to be the $L^p+L^\infty$ norm of
$$
(1+|x|^2)^{\alpha/2}V(x)~.
$$
Let
$$
a(\zeta)=VR(\zeta)~.
$$
In the following, $C$ denotes a numerical constant. Choose $s$ such that $1<s<\alpha/2$.

If $\eta={\rm Im}\, \zeta \neq 0$, then as an operator on $L^2(\mathbb{R}^3)$, $a(\zeta)$ is analytic in $\zeta$, and
\begin{equation}\label{A.1}
||a(\zeta)||\leq C\beta(V)\eta ^{-\frac{1}{4}}
\end{equation}
for $\eta \geq 1$, say.
If $a(\zeta)$ is regarded as an operator on the weighted space $L_{s}^{2}(\mathbb{R}^3)$, for $1<s<\alpha/2$, then $a(\zeta)$ is analytic for Im $\zeta\neq 0$, and
\begin{equation}\label{A.2}
||a(\zeta)||_s \leq C\beta(V)
\end{equation}
and
\begin{equation}\label{A.3}
||a(\zeta')-a(\zeta||_s \leq C\beta(V)|\zeta-\zeta'|^\gamma
\end{equation}
uniformly in $\zeta$ and $\zeta'$, where $\gamma > 0$. For fixed $\epsilon \geq 0$,
\begin{equation}\label{A.4}
\lim_{\lambda\rightarrow \pm\infty}||a(\lambda\pm i\epsilon)||_s=0~.
\end{equation}
Finally, $a(\zeta)$ is {\em compact} on $L_{s}^{2}$. [4, Proposition 3.1]

Defining $f_0(\lambda; x)=\langle\delta(H-\lambda)x,x\rangle$ for $x\in L_{s}^{2}$ and real, we have
\begin{equation}\label{A.5}
f_0(\lambda,x)\leq C||x||_{s}^{2}
\end{equation}
uniformly in $\lambda$. Thus, $f_0$ is a spectral form for $H$. [4, Proposition 2.1]

If $x\in L_{s}^{2}(\mathbb{R}^3)$ and $f_0(\lambda,x)=0$, then $\psi=(H-\lambda)^{-1}x\in L_{s}^{2}(\mathbb{R}^3)$ and
\begin{equation}\label{A.6}
||\psi||\leq C(\zeta)||x||_s~.
\end{equation}
where $\lambda\geq \delta >0$. [4, Proposition 2.6]

For $x\in L^{2}(\mathbb{R})$ and Im$\zeta\neq 0$, $\psi=R(\zeta)x\in L_{s}^{2}(\mathbb{R}^3)$ and
\begin{equation}\label{A.7}
||\psi||_s\leq C(\zeta)||x||~.
\end{equation}
We have given no reference for (\ref{A.1}) or (\ref{A.7}). The former can be proved by estimating the Hilbert-Schmidt norm of $a(k^2)$, using that $R(k^2)$ is convolution by
$$
(4\pi|x|)^{-1}e^{ik|x|}
$$
for Im$\,k>0$. The latter follows by noting that $WR(\zeta)$ is bounded on $L^{2}(\mathbb{R}^3)$ for $W\in L^2+L^{\infty}$.

%


\end{document}